
\input amstex
\documentstyle{amsppt}
\def\const{\operatorname{const}}
\rightheadtext {On the separate algebraicity \dots}
\topmatter
\title
On the separate algebraicity along the families of algebraic curves.
\endtitle
\author
R.A.~Sharipov and E.N.~Tsyganov
\endauthor
\abstract
     A new generalization of the classical separate algebraicity
theorem is suggested and proved.\par
\endabstract
\address
Department of Mathematics, Bashkir State University,
Frunze street 32, 450074 Ufa, Russia
\endaddress
\email
yavdat\@bgua.bashkiria.su
\endemail
\thanks
Paper is written under the financial support of International
Scientific Foundation of Soros, projects \#RK4000 and \#RK4300,
and Russian Fund RFFI, project \#93-01-00273.
\endthanks
\endtopmatter
\document
\head
1. Introduction.
\endhead
     The theorem on separate holomorphy (Hartogs theorem) is
one of the basic theorems in the theory of functions of several
complex variables (see, for example, in \cite{1}). It's formulated as
follows.
\proclaim {Theorem 1} Let $f(z)=f(z_1,\ldots,z_n)$ be a function
in the domain $D\subset\Bbb C^n$, which is holomorphic in each
variable $z_i$ for any fixed values of other variables. Then it is
holomorphic function in $D$.
\endproclaim
     Theorem~1 has certain modifications for polynomial,
rational and algebraic functions (see \cite{2}). In the last case
we have the following statement.
\proclaim {Theorem 2} Let $f(z)=f(z_1,\ldots,z_n)$ be a holomorphic
function in the domain $D\subset\Bbb C^n$, which is algebraic in
each variable $z_i$ for any fixed values of other variables. Then
it is holomorphic branch in $D$ for some algebraic function
given by a polynomial equation $P(f,z_1,\ldots,z_n)=0$.
\endproclaim
     Let's study the assumptions of these two theorems. In both
cases we may consider the following families of complex lines,
one per each variable:
$$
\cases
z_i=c^{(m)}_i=\const\quad&\text{for\ }i\neq m,\\
z_m=t\in\Bbb C &\text{for\ }i=m.
\endcases
\tag1.1
$$
These are the coordinate lines represented in a parametric form,
$t\in\Bbb C$ is a parameter. Restricting $f(z)$ to such lines,
we obtain the following functions:
$$
f_m(t)=f(c^{(m)}_1,\ldots,c^{(m)}_{m-1},t,c^{(m)}_{m+1},\ldots
,c^{(m)}_n),
\tag1.2
$$
which are holomorphic in $t$ under the assumptions of theorem~1
and algebraic in $t$ under the assumptions of theorem~2.\par
     In terms of coordinate lines \thetag{1.1} the above classical
theorem~2 can be formulated as follows: every holomorphic function
f(z), which is algebraic along each coordinate line in $D$, is
algebraic as a function of several complex variables.
In the present paper we consider some generalized versions
of this statement, when coordinate lines are replaced by certain
classes of complex curves. The necessity of such generalization was
fist recognized by A.B.~Sukhov \cite{3, 4} in connection with
the problem of Poincar\'e--Alexander (see \cite{5--8}) and
its generalizations (see \cite{9--13}).
      As a first step in \cite{3, 4} coordinate lines were replaced
by quadrics. Second step was done in \cite{14}. Here instead of
coordinate lines arbitrary algebraic curves are considered. They
are taken in a parametric form with parameter $t\in\Bbb C$:
$$
\left\{
\aligned
&z_1=R^{(m)}_1 (t,c^{(m)}_1,\ldots,c^{(m)}_{n-1}),\\
&\text{\hbox to 4.5cm{\leaders\hbox to 1em{\hss.\hss}\hfill}}\\
&z_n=R^{(m)}_n(t,c^{(m)}_1,\ldots,c^{(m)}_{ n-1}).
\endaligned\right.
\tag1.3
$$
Algebraic curves \thetag{1.3} form $n$ families $R^{(1)},\ldots,
R^{(n)}$, index $m$ numerates these families. Coordinate functions
$R^{(m)}_i$ in \thetag{1.3} are algebraic in $t$. Other complex
parameters $c^{(m)}_1,\ldots,c^{(m)}_{ n-1}$ are introduced in
order to determine particular curve within $m$--th family.\par
     Note that in general the whole set of curves \thetag{1.3}
can't be mapped to the coordinate lines by an algebraic map.
This is why we should develop some special tools for this case.
\par
     Given a function $f(z)$ in $D$, we may consider its restriction
to the curves \thetag{1.3}. This gives us $n$ functions analogous to
\thetag{1.2}:
$$
f_m(t)=f(R^{(m)}_1(t,c^{(m)}_1,\ldots,c^{(m)}_{ n-1}),\ldots,
R^{(m)}_n(t,c^{(m)}_1,\dots,c^{(m)}_{ n-1}))
\tag1.4
$$
We say that the function $f(z)$ is {\it algebraic along the curves}
\thetag{1.3}, if each function $f_m(t)$ in \thetag{1.4} depends
algebraically on $t$. In terms of curves \thetag{1.3} the
separate algebraicity theorem can be stated as follows.
\proclaim {Theorem 3} Let $f(z)=f(z_1,\ldots,z_n)$ be a holomorphic
function in a domain $D\subset\Bbb C^n$, which is algebraic
along algebraic curves forming $n$ regular\footnotemark"$^{1)}$"\
families being in general position\footnotemark"$^{2)}$".
Then it is holomorphic branch in $D$ for some algebraic function
given by a polynomial equation $P(f,z_1,\ldots,z_n)=0$.
\endproclaim
\footnotetext"$^{1,2)}$"{\ Precise definitions see in section~2.}
     This theorem~3 was proved in \cite{14}, but in a more
restricted form. Coordinate functions $R^{(m)}_i$ in
\thetag{1.3} were assumed to be algebraic not only in $t$,
but in whole set of their arguments $t,c^{(m)}_1,\dots,
c^{(m)}_{ n-1}$. In the present paper we exclude this extra
hypothesis concerning $c^{(m)}_1,\dots,c^{(m)}_{ n-1}$,
but make more restrictive assumptions concerning the main
variable $t$ in \thetag{1.3} and \thetag{1.4}. The functions
$z_1(t),\ldots,z_n(t)$ and $f_m(t)$ are supposed to be
polynomials in $t$. Below is the main theorem that we prove in
this paper.
\proclaim {Theorem 4} Let $f(z)=f(z_1,\ldots,z_n)$ be a holomorphic
function in a domain $D\subset\Bbb C^n$, which is polynomial
along polynomial curves forming $n$ regular families being
in general position. Then it is holomorphic branch in $D$ for
some algebraic function given by a polynomial equation
$P(f,z_1,\ldots,z_n)=0$.
\endproclaim
     As for the most general theorem~3 it still remains an
unproved conjecture. The proof of this theorem is the
subject for the next publication.\par
     Our paper is organized as follows. In section~2 we give
precise definitions and preliminary discussions. Sections 3
contains the sketch of the proof of the main theorem~4. Here we
omit all details and technicalities. Other sections are devoted to
the strict proofs for all auxiliary results we need to prove the
main theorem.\par
\head
2. Precise definitions and some preliminaries.
\endhead
     Parametric equations \thetag{1.3} define our basic families
of polynomial curves. To fit the theorem~4 they should be regular
in the sense of the following definition.
\definition{Definition 1} We say that $m$--th family of polynomial
curves \thetag{1.3} is regular in a domain $D\subset\Bbb C^n$,
if corresponding map
$$
R^{(m)}:(t,c^{(m)}_1,\dots,c^{(m)}_{ n-1})\to (z_1,\ldots,z_n)
\in D\tag2.1
$$
in \thetag{1.3} is the biholomorphic diffeomorphism of some
domain $U_m\subset\Bbb C^n$ and $D$.
\enddefinition
If $m$--th family is regular, we may treat $t,c^{(m)}_1,\dots,
c^{(m)}_{ n-1}$ as new curvilinear complex coordinates in $D$.
Then the following map
$$
\varphi^{(m)}_\tau:(t,c^{(m)}_1,\dots,c^{(m)}_{ n-1})\to
(t+\tau,c^{(m)}_1,\dots,c^{(m)}_{ n-1})
\tag2.2
$$
is a translation along the curves of $m$--th family and $\tau$
is a magnitude of translation. By means of diffeomorphism
\thetag{2.1} we may rewrite \thetag{2.2} in cartesian
coordinates. Here this map $\varphi^{(m)}_\tau: z\to\tilde z$
is given by holomorphic functions
$$
\left\{
\aligned
&\tilde z_1=\varphi^{(m)}_1(\tau,z_1,\ldots,z_n),\\
&\text{\hbox to 4cm{\leaders\hbox to 1em{\hss.\hss}\hfill}}\\
&\tilde z_n=\varphi^{(m)}_n(\tau,z_1,\ldots,z_n).
\endaligned\right.
\tag2.3
$$
{}From \thetag{2.2} we find that $\varphi^{(m)}_\tau(
\varphi^{(m)}_\theta(z))=\varphi^{(m)}_{\tau+\theta}(z)$.
Therefore $\varphi^{(m)}_\tau$ define the holomorphic
one--parametric local group of diffeomorphisms in $D$.
This local group generates the holomorphic vector field
$\bold X_m$ in $D$, tangent to the curves of $m$--th
family. Under the assumptions of theorem~4 we have
$n$ holomorphic vector--fields:
$$
\bold X_1,\ldots,\bold X_n.\tag2.4
$$
      In terms of curvilinear coordinates $t,c^{(m)}_1,
\ldots,c^{(m)}_{ n-1}$ in $\tilde D$ the $m$-th vector--field
$\bold X_m$ from \thetag{2.4} is represented by the following
differential operator:
$$
\bold X_m=\frac{\partial}{\partial t}.\tag2.5
$$
\definition{Definition 2} Given $n$ regular families of polynomial
curves in $D$ we shall say that they are in general position, if
corresponding vector--fields \thetag{2.4} are linearly independent
at each point $z\in D$.
\enddefinition
     Now let's consider the polynomial curves \thetag{1.3}, which
are supposed to form $n$ regular families in general position.
Denote by $d_{1m},\ldots,d_{nm}$ and $p_m$ the degrees of
polynomials $z_1(t),\ldots,z_n(t)$ and $f_m(t)$ in \thetag{1.3}
and \thetag{1.4}. They depend on $c^{(m)}_1,\ldots,c^{(m)}_{n-1}$.
Because of regularity of $m$--th family of curves in $D$
one can treat them as  functions of a point $z\in D$:
$$
\aligned
d_{im}(z)&=\deg R^{(m)}_i(t),\ i=1,\ldots,n,\ m=1,\ldots,n,\\
p_m(z)&=\deg f_m(t),\ m=1,\ldots,n.
\endaligned
\tag2.6
$$
Functions \thetag{2.6} are integer--valued functions in $D$,
which remain constant along the curves of appropriate family.
In section~4 we shall prove the following theorem.
\proclaim{Theorem 5} Under the assumptions of theorem~4 one can
find a smaller subdomain $\tilde D\subset D$, such
that all functions \thetag{2.6} are constant in $\tilde D$.
\endproclaim
\noindent Now we are ready to explain the idea for the proof of
the main theorem~4.
\head
3. Proof of the main theorem~4.
\endhead
     Note that in theorems~5 we contract the domain $D$ to
$\tilde D\subset D$. This doesn't affect the ultimate result,
since holomorphic function $f(z)$ in $D$, which is algebraic
in some smaller subdomain of $D$, is algebraic in $D$. Therefore,
we can reduce our domain as many times as we need.\par
     According to the theorem~6 we have the subdomain $\tilde D
\subset D$, where the degrees of polynomials in \thetag{1.3} and
\thetag{1.4} are constant:
$$
\xalignat 2
&\deg R^{(m)}_i(t)=d_{im},
&
&\deg f_m(t)=p_m.
\endxalignat
$$
Choose some other non--negative integer numbers $k_{1},\ldots,k_{n}$
and $q$ and consider the following monomial in $\tilde D$:
$$
M(k_1,\ldots,k_n,q)=f^q\cdot(z_1)^{k_1}\cdot\ldots\cdot(z_n)^{k_n}.
\tag3.1
$$
where $f=f(z)$. Considered as a function in $\tilde D$ monomial
\thetag{3.1} is holomorphic. We restrict it to the some arbitrary
curve of $m$--th family in $\tilde D$. Then $M(k_1,\ldots,k_n,q)$
becomes a polynomial in $t$, its degree is given by the formula:
$$
q_m=\deg M(k_1,\ldots,k_n,q)=q\,p_m+\sum^n_{i=1}k_i\,d_{im}.
\tag3.2
$$
If $N>q_m$ is some integer number, then, applying $N$-th power
of the differential operator \thetag{2.5} to the monomial
\thetag{3.1}, we get
$$
(\bold X_m)^N M(k_1,\ldots,k_n,q)=0.
\tag3.3
$$
Respective to $M(k_1,\ldots,k_n,q)$ the equality \thetag{3.3}
is a linear differential equation of $N$--th order with
holomorphic coefficients. Now we unite all these equations
\thetag{3.3} into the system
$$
(\bold X_m)^N\varphi(z)=0,\ m=1,\ldots,n,
\tag3.4
$$
and denote by $V(N)$ the set of their holomorphic solutions
$\varphi(z)$ in $\tilde D$:
$$
V(N)=\{\varphi\in\Cal H(\tilde D): (\bold X_m)^N\varphi=0
\text{\ for all\ } m=1,\ldots,n\}.\tag3.5
$$
Since the equations \thetag{3.4} are linear and homogeneous, the
set of their solutions $V(N)$ forms a linear space over the field
of complex numbers.
\proclaim{Theorem 6} For any integer $N>0$ the complex linear
space \thetag{3.5} has finite dimension and $\dim V(N)\leq N^n$.
\endproclaim
In section~4 we shall prove this theorem giving the
estimate for $\dim V(N)$. Now denote by $M(N)$ the number of
monomials \thetag{3.1}, for which the degrees \thetag{3.2} of
their restrictions to the curves are less than $N$, i.e.
$q_m<N$ for all $m=1,\ldots,n$. If for some value of $N$ we
find that $M(N)>N^n$, then we obtain that certain number
of monomials \thetag{3.1} are linearly dependent over $\Bbb C$.
This will give a polynomial equation $P(f,z^1,\ldots,z_n)=0$
for the function $f(z)$ and will terminate the proof of theorem~4.
\par
     According to the above conclusion the last step in the proof
of theorem~4 is to be the proper estimate for $M(N)$ at least
for some certain value of $N$. Suppose $p=\max\{p_1,\ldots,p_n\}$
and $d_i=\max\{d_{i1},\ldots,d_{in}\}$ (see formula \thetag{3.2}).
Choose some arbitrary integer number $K>1$ and let
$$
N=N(K)=1+K\,p+\sum^n_{i=1}K\,d_i.
\tag3.6
$$
Then for $q=1,\ldots,K$ and for $k_i=1,\ldots,K$ the degrees
of corresponding monomials \thetag{3.1} are less than $N$.
For the number of such monomials we have
$$
M(N)\geq K^{n+1}.
$$
For the dimension of V(N) from \thetag{3.6} we derive another
estimate
$$
\dim V(N)\leq N^n\leq\const\cdot K^n,\text{\ as\ }K\to\infty.
$$
Comparing these to estimates we conclude that $M(N)>\dim V(N)$
for some large enough value of $K$.\par
     So the main theorem~4 is proved provided the
theorems~5~and~6 hold. The rest part of paper is devoted to proof
of these two theorems.\par
\head
4. Proof of the theorems 5 and 6.
\endhead
     Let's study the integer-valued functions $d_{im}(z)$ and
$p_m(z)$ in \thetag{2.6}. The whole set of these functions
can be treated as a map
$$
\nu: D\to\Bbb Z^r,\tag 4.1
$$
where $r=n(n+1)$. Here in \thetag{4.1} $S=\Bbb Z^r$ is a countable
set. Therefore we may apply the following theorem to the map
\thetag{4.1}.
\proclaim{Theorem 7} For any map $\nu: D\to S$ from some domain
$V\subset\Bbb C^{n-1}$ to the countable set $S$ one can find
at least one value $s\in S$, such that its inverse image
$\nu^{-1}(s)=\{v\in V:\ \nu(v)=s\}$ is a dense set in some
subdomain $\tilde D\subset D$.
\endproclaim
     The domain $\tilde D$ given by the theorem 7 is exactly
the domain we need to prove the theorem 5. So we reduced theorem 5
to the topological theorem 7. We shall not prove this topological
theorem, since it is direct consequence of the theorem of Baire
(see \cite{15} chapter 3 section 5 theorem 8 and \cite{16}
chapter 2 section 3). Similar theorems were used in \cite{2} and
\cite{14}.\par
     The next is the proof of the theorem 6. Let's consider the
vector fields \thetag{2.4} and their local one--parametric
groups \thetag{2.3}. Fix some point $z_0\in D$ and apply
the first map \thetag{2.2} to this point. As a result we obtain
some one--parametric set of points $z(t_1)=\varphi^{(1)}_{t_1}(z_0)$,
i.e. the holomorphic complex curve in $D$. Then we apply another
map $\varphi^{(2)}_{t_2}$ to the points of this line and obtain
a complex surface
$$
z(t_1,t_2)=\varphi^{(2)}_{t_2}(z(t_1)).
$$
Continuing this process on the $n$--th step we obtain the map
$$
(t_1,\ldots,t_n)\to z=z(t_1,\ldots,t_n)\in D\tag4.2
$$
Since curves \thetag{1.3} are in general position, this map
\thetag{4.2} can be treated as holomorphic change of coordinates
in some neighborhood $\tilde D$ of the point $z_0=z(0,\ldots,0)$.
\par
     In general case vector fields $\thetag{2.4}$ do not
commutate. Therefore their local groups of diffeomorphisms
\thetag{2.2} also don't commute. Nevertheless we can
use the map \thetag{4.2} to simplify the differential equations
\thetag{3.4}. Because of \thetag{2.5} the last $n$--th equation
\thetag{3.5} is written as:
$$
\left(\frac{\partial}{\partial t_n}\right)^N
\varphi(t_1,\ldots,t_{n-1},t_n)=0.
$$
The previous $(n-1)$--th equation \thetag{3.5} can be simplified
only for $t_n=0$:
$$
\left(\frac{\partial}{\partial t_{n-1}}\right)^N
\varphi(t_1,\ldots,t_{n-1},0)=0.
$$
Continuing this process backward from $n$--th equation to the first
equation \thetag{3.5}, we receive the following system of differential
equations for the function $\varphi(z)$:
$$
\left\{
\aligned
&\left(\frac{\partial}{\partial t_n}\right)^N
\varphi(t_1,t_2,\ldots,t_{n-1},t_n)=0,\\
&\text{\hbox to 5cm{\leaders\hbox to 1em{\hss.\hss}\hfill}}\\
&\left(\frac{\partial}{\partial t_1}\right)^N
\varphi(t_1,0,\ldots,0,0)=0.
\endaligned
\right.
\tag4.3
$$
     One can easily derive the formula for the general solution of
the system \thetag{4.3}. It is the following polynomial with
arbitrary constant coefficients:
$$
\varphi=\sum^{N-1}_{i_1=0}\ldots\sum^{N-1}_{i_n=0}
C_{i_1,\ldots,i_n}\,(t_1)^{i_1}\cdot\ldots\cdot(t_n)^{i_n}
\tag4.4
$$\par
     From \thetag{4.4} one can easily count the number of
linearly independent solutions for \thetag{4.3}. It's exactly
$N^n$. But the number of linearly independent solutions of
\thetag{3.5} can be less than $N^n$. The equations \thetag{4.3}
follow from \thetag{3.5}, but they are not completely equivalent
to \thetag{3.5}, except for the case when vector--fields
\thetag{2.3} are commuting. The estimate $\dim V(N)\leq N^n$
is proved. This completes the proof of the theorem 6.
\head
5. Acknowledgements.
\endhead
     Authors are very grateful to S.I.~Pinchuk and A.B.~Sukhov who
keep hot our interest for the discussed problems.

\enddocument
\end